\documentclass[a4paper, 11pt]{article}
\usepackage{amsmath, amssymb, eucal, amscd, amstext}

\newtheorem{thm}{Theorem}[section]
\newtheorem{lem}[thm]{Lemma}
\newtheorem{eg}[thm]{Example}
\newtheorem{prop}[thm]{Proposition}
\newtheorem{cor}[thm]{Corollary}
\newtheorem{defn}[thm]{Definition}
\newtheorem{rem}[thm]{Remark}
\newtheorem{ntn}[thm]{Notation}

\newenvironment{prf}{{\noindent \textbf{Proof:}\ }}{\hfill $\Box$\\ \smallskip}

\numberwithin{equation}{section}

\newcommand{\abs}[1]{{\lvert#1\rvert}}

\newcommand{\id}{{\rm id}}

\newcommand{\la}{\langle}
\newcommand{\ra}{\rangle}
\newcommand{\smnoind}{\smallskip\noindent}
\newcommand{\CL}{\mathcal{L}}

\newcommand{\sbm}{{\rm Bimod}^*}
\newcommand{\sph}{\mathfrak{S}_1}
\newcommand{\F}{\mathfrak{F}}
\newcommand{\RP}{\mathbb{R}_+}
\newcommand{\FH}{\mathcal{H}}
\newcommand{\FK}{\mathcal{K}}
\newcommand{\CC}{\mathcal{C}}

\begin{document}

\title{Property $(T)$ and strong Property $(T)$ for unital $C^*$-algebras
\footnote{This work is jointly supported by
Hong Kong RGC Research Grant (2160255), Hong Kong RGC Direct Grant (2060319), the National Natural Science Foundation of China (10771106) and NCET-05-0219}}
\author{Chi-Wai Leung and Chi-Keung Ng}

\maketitle
\begin{abstract}
 
In this paper, we will give a thorough study of the notion of Property $(T)$ for $C^*$-algebras (as introduced by M.B. Bekka in
\cite{Bek-T}) as well as a slight stronger version of it, called ``strong property $(T)$'' (which is also an analogue of the corresponding concept in the case of discrete groups and type $\rm II_1$-factors).
More precisely, we will give some interesting equivalent formulations as well as some permanence
properties for both property $(T)$ and strong property $(T)$. 
We will also relate them to certain $(T)$-type properties of the unitary group of the underlying $C^*$-algebra.

\medskip

\end{abstract}

\medskip

\medskip

\section{Introduction}

\medskip

Property $(T)$ for locally compact groups was first defined by D. Kazhdan in \cite{Kaz} and was later extended to Hausdorff topological groups. 
In \cite{Mar}, Property $(T)$ for a pair of groups $H\subseteq G$ was introduced. 
This notion was proved to be very useful and was studied by many people (see e.g. \cite{Bek-unit}, \cite{BV}, \cite{HV}, \cite{Jol}, \cite{Kaz}, \cite{Mar} and \cite{V}). 

\medskip

In \cite{Con}, A. Connes introduced the notion of property $(T)$ for type $\rm II_1$-factors and this notion was then extended to von Neumann algebras in \cite{CJ}. 
A discrete group $G$ has property $(T)$ if and only if the von Neumann algebra generated by the left regular representation of $G$ has property $(T)$ (this was first proved in \cite{CJ} for ICC groups and was generalized by P. Jolissaint in \cite{Jol-disc} to general discrete groups). 
The notion of property $(T)$ for a pair of von Neumann algebras was defined by S. Popa in \cite{Popa}. 
This notion was also proved to be very useful in the study of von Neumann algebras. 

\medskip

Recently, M.B. Bekka introduced in \cite{Bek-T} the interesting notion of property $(T)$ for a pair consisting of a unital $C^*$-algebra and a unital $C^*$-subalgebra. 
He showed that a countable discrete group $G$ has property $(T)$ if and only if its full (or equivalently reduced) group $C^*$-algebra has property $(T)$. 
In \cite{Bro}, N.P. Brown did a study of property $(T)$ for $C^*$-algebras and showed that a nuclear unital $C^*$-algebra $A$ has property $(T)$ if and only if $A = B\oplus C$ where $B$ is finite dimensional and $C$ admits no tracial state. 

\medskip

The aim of this paper is to give a thorough study of property $(T)$ as well as a slightly stronger version called strong property $(T)$ for unital $C^*$-algebras. 
On our way, we will show that our stronger version is equally good (if not a better) candidate for the notion of property $(T)$ for a pair of unital $C^*$-algebras. 

\medskip

The paper is organised as follows. 
In Section 2, we will give two simple and useful reformulations of both property $(T)$ and strong property $(T)$. 
In section 3, we consider two Kazhdan constants $t^A_u$ and $t^A_c$ for a $C^*$-algebra $A$ which are the analogous of the Kazhdan constant for
locally compact groups (see \cite{V}). 
We will show that $A$ has property $(T)$ (respectively, strong property $(T)$) if and only if $t_c^A > 0$ (respectively, $t_u^A > 0$). 
Through them, we obtain some interesting reformulations of property $(T)$ and strong property $(T)$. 
In particular, we show that one can check property $(T)$ by looking at just one bimodule. 
We will also show that one can express the gap between property $(T)$ and strong property $(T)$ by another Kazhdan constant $t^A_s$. 

\medskip

In section 4, we obtain some permanence properties for
property $(T)$ and strong property $(T)$, including quotients, direct sums, tensor products and crossed products.
In section 5, we will show that finite dimensional $C^*$-algebras have strong property $(T)$. 
Moreover, we show that a corresponding result of Bekka concerning relation between property $(T)$ of discrete groups and their group $C^*$-algebras as well as a corresponding result of Brown concerning amenable property $(T)$ $C^*$-algebras also holds for strong property $(T)$. 
In Section 6, we study the relation between property $(T)$ (as well as strong property $(T)$) of a unital $C^*$-algebra $A$ and certain $(T)$-type properties of the unitary group of $A$. 

\medskip

Let us first set the following notations that will be used throughout the whole paper. 

\medskip

\begin{ntn}
(1) $A$ is a unital $C^*$-algebra and $B \subseteq A$ is a $C^*$-subalgebra containing the identity of $A$. 
Set $A^{Dou} := A\otimes_{\rm max} A^{\rm op}$ (where $A^{\rm op}$ is the ``opposite $C^*$-algebra'' with $a^{\rm op}b^{\rm op} = (ba)^{\rm op}$).

\smnoind
(2) $\F(E)$ is the set of all non-empty finite subsets of a set $E$ and $\sph(X)$ is the unit sphere of a normed space $X$.

\smnoind 
(3) $U(A)$ and $S(A)$ are respectively the unitary group and the state space of
$A$. 

\smnoind 
(4) $\sbm(A)$ is the collection of unitary
equivalence classes of unital Hilbert bimodules over $A$ (or
equivalently, non-degenerate representations of $A^{Dou}$). 
For any $H\in \sbm(A)$, let
$$H^B\ :=\ \{\xi\in H: b\cdot \xi = \xi \cdot b\ {\rm for\ all}\ b\in B\}$$
and $P^B_H: H \rightarrow H^B$ be the orthogonal projection.
Elements in $H^B$ are called \emph{central vectors for $B$}. 
Moreover, for any $(Q,\beta) \in \F(A)\times \RP$, set
$$V_H(Q,\beta)\ :=\ \{\xi\in \sph(H): \|x\cdot \xi - \xi\cdot x\|
< \beta\ {\rm for\ all}\ x\in Q\}.$$
Elements in $V_H(Q,\beta)$ are called \emph{$(Q,\beta)$-central unit vectors}. 
On the other hand, a net of vectors
$(\xi_i)_{i\in I}$ in $\sph(H)$ is called an \emph{almost central unit vector for $A$} if $\|a\cdot \xi_i - \xi_i \cdot a\|
\rightarrow 0$ for any $a\in A$.

\smnoind 
(5) For any topological group $G$, we denote by ${\rm
Rep}(G)$ the collection of all unitary equivalence classes of continuous
unitary representations of $G$. If $(\pi, H)\in {\rm Rep}(G)$, we
let
$$H^G\ :=\ \{\xi\in H: \pi(s) \xi = \xi\ {\rm for\ all}\ s\in G\}$$
and $P^G_H: H \rightarrow H^G$ be the orthogonal projection.
Furthermore, if $F\in \F(G)$ and $\epsilon > 0$, we set
$$V_\pi(F,\epsilon) = \{ \xi\in \sph(H): \|\pi(t)\xi -\xi\|<
\epsilon\ {\rm for\ all}\ t\in F\}.$$

\smnoind 
(6) For any $(\mu,H), (\nu,K) \in {\rm Rep}(G)$, we write
$(\mu, H) \leq (\nu, K)$ if $(\mu,H)$ is a subrepresentation of
$(\nu,K)$.
\end{ntn}

\medskip

\medskip

\section{Definitions and basic properties}

\medskip

Let us first recall Bekka's notion of property $(T)$ in \cite{Bek-T}. 
The pair $(A, B)$ is said to have \emph{property $(T)$} if there
exist $F\in \F(A)$ and $\epsilon > 0$ such that for any $H\in \sbm(A)$, if $V_H(F, \epsilon)\neq \emptyset$, then $H^B\neq (0)$. 
In this case, $(F, \epsilon)$ is called a \emph{Kazhdan pair} for $(A, B)$. 
Moreover, $A$ is said to have \emph{property $(T)$} if the pair $(A,A)$ has property $(T)$. 

\medskip

Note that Bekka's definition comes from the original definition of property $(T)$ for groups (see e.g. \cite[Definition 1.1(1)]{Jol}). 
We will now give a slightly stronger version which comes from an equivalent form of property $(T)$ for groups (see \cite[Theorem 1.2(b2)]{Jol}). 
Note that the corresponding stronger version of property $(T)$ for type $\rm II_{1}$-factor is also equivalent to property $(T)$ (see e.g. \cite[Proposition 1]{CJ}) but we do not know if it is the case for $C^*$-algebras. 

\medskip

\begin{defn}
\label{def-sT} 
The pair $(A,B)$ is said to have \emph{strong property $(T)$} if for any $\alpha > 0$, there exist $Q\in \F(A)$ and $\beta > 0$ such that for any $H\in \sbm(A)$ and any $\xi\in V_H(Q, \beta)$, we have $\left\|\xi - P^B_H(\xi)\right\| < \alpha$. 
In this case, $(Q,\beta)$ is called a \emph{strong Kazhdan pair} for $(A, B, \alpha)$. 
We say that $A$ has strong property $(T)$ if $(A,A)$ has such property.
\end{defn}

\medskip

It is clear that if $A$ has property $(T)$ (respectively, strong property $(T)$) then so has the pair $(A,B)$. 
Moreover, by taking $\alpha < 1/2$, we see that strong property $(T)$
implies property $(T)$. 
We will see later that strong property $(T)$ is an equally good (if not a better) candidate for the notion of property $(T)$ for a pair of $C^*$-algebras. 

\medskip

Let us now give the following simple reformulation of property $(T)$ and strong property $(T)$ which will be useful in Section 6.

\medskip

\begin{lem}
\label{F(U)} For any $(Q,\beta)\in \F(A)\times \RP$, there exists
$(Q',\beta')\in \F(U(A))\times \RP$ such that $V_H(Q',\beta')
\subseteq V_H(Q,\beta)$ for any $H\in \sbm(A)$. Consequently, one
can replace $\F(A)$ by $\F(U(A))$ in the definitions of both property $(T)$ and strong property $(T)$.
\end{lem}
\begin{prf}
This lemma is clear if $Q = \{0\}$. Let $Q\setminus \{0\} =
\{x_1,...,x_n\}$ and $M = \max\{\|x_1\|,...,\|x_n\|\}$. For each
$k\in \{1,...,n\}$, consider $u_k,v_k\in U(A)$ such that $2x_k =
\|x_k\|((u_k + u_k^*) + i(v_k + v_k^*))$. If we take $Q'$ to be the
set $\{u_1,u_1^*,v_1,v_1^*,...,u_n,u_n^*,v_n,v_n^*\}$ and $\beta' =
\frac{\beta}{2M}$, then $V_H(Q',\beta')\subseteq V_H(Q,\beta)$ for
any $H\in \sbm(A)$.
\end{prf}

\medskip

Before we give a second simple reformulation, we need to set some notations. 
Let $S(D)$ and $S_t(D)$ be respectively the sets of all states and
the set of all tracial states on a $C^*$-algebra $D$. For any
$\tau\in S(D)$ and any cardinal $\alpha$, we denote by $M_\tau$ the
GNS construction for $\tau$ and by $M_{\tau, \alpha}$ the
$\alpha$-times direct sum, $\bigoplus_\alpha M_\tau$, of $M_\tau$
(we use the convention that $\bigoplus_0 M_\tau = \{0\}$).

\medskip

\begin{defn}
\label{def-tq} 
Let
$$\FH\ :=\ \bigoplus_{\tau\in S(A^{Dou})} M_\tau \qquad {\rm and}
\qquad \FK\ :=\ \bigoplus_{\tau\in S_t(A)} M_\tau.$$
We called $\FH$ and $\FK$ the \emph{universal} and the
\emph{standard} bimodules (over $A$) respectively. Moreover, a
bimodule of the form $\bigoplus_{\tau \in S_t(A)} M_{\tau,
\alpha_\tau}$ is called  a \emph{quasi-standard bimodule}.
\end{defn}

\medskip

\begin{prop}
\label{acv} (a) $(A,B)$ has property $(T)$ if and only if for any
$H\in \sbm(A)$, the existence of an almost central unit
vector for $A$ in $H$ will imply that $H^B\neq \{0\}$.

\smnoind (b) The following statement are equivalent.
\begin{enumerate}
\item[\rm (i)] $(A,B)$ has strong property $(T)$.
\item[\rm (ii)] For any almost central unit vector $(\xi_i)_{i\in I}$ for $A$ in any bimodule $H\in \sbm(A)$, we have 
$\left\|\xi_i - P^B_H(\xi_i)\right\| \rightarrow 0$.
\item[\rm (iii)] For any almost central unit vector $(\xi_i)_{i\in I}$ for $A$ in $\FH$ and any $n\in \mathbb{N}$, there exists $i_n\in I$ with $\left\|\xi_{i_n} - P^B_{\FH}(\xi_{i_n})\right\| < 1/n$.
\end{enumerate}
\end{prop}
\begin{prf}
(a) This part is well known.

\smnoind (b) It is clear that (i)$\Rightarrow$(ii) and (ii)$\Rightarrow$(iii).
To obtain (iii)$\Rightarrow$(i), we suppose, on the contrary, that $(A,B)$ does not have strong property $(T)$. 
Then one can find $\alpha_0 > 0$ such that for any $i = (Q,\beta) \in I :=
\F(A)\times \RP$, there exist $H_i\in \sbm(A)$ and $\xi_i\in V_{H_i}(Q,\beta)$ with $\|\xi_i - P^B_{H_i}(\xi_i)\| \geq \alpha_0$. 
If $K_i = \overline{A\cdot \xi_i \cdot A}$, then $H_i = K_i \oplus K_i^\bot$ and $H_i^B = K_i^B \oplus (K_i^\bot)^B$. As $\xi_i\in K_i$, we have 
$$\|\xi_i - P^B_{K_i}(\xi_i)\|\ =\ \|\xi_i - P^B_{H_i}(\xi_i)\|\ \geq\ \alpha_0.$$ 
We set $\mathcal{X} := \{K_i: i\in I\} \subseteq \sbm(A)$ and $K_0 := \bigoplus_{K\in \mathcal{X}} K$. 
Since all bimodules in $\mathcal{X}$ are cyclic (as representations of $A^{Dou}$) and any two elements in $\mathcal{X}$ are inequivalent, $K_0$ is a Hilbert sub-bimodule of $\FH$. 
Moreover, each $K_i$ is equivalent to a unique element in $\mathcal{X}$ and this gives a canonical Hilbert
bimodule embedding $\Psi_i: K_i \rightarrow K_0$. It is easy to check that
$(\Psi_i(\xi_i))_{i\in I}$ is an almost central unit vector for $A$ in $\FH$ with
$\|\Psi_i(\xi_i) - P^B_{\FH}(\Psi_i(\xi_i))\| = \|\xi_i - P^B_{K_i}(\xi_i)\| \geq \alpha_0$
for every $i\in I$ (since $\Psi_i(K_i)$ is a direct summand of $\FH$). This contradicts Statement (iii).
\end{prf}

\medskip

Since Proposition \ref{acv} is so fundamental to our discussions, we may use it without mentioning it explicitly throughout the whole paper. 

\medskip

\begin{rem}
\label{rem-acv} 
(a) Note that if $A$ is separable, then in the above proposition, one can replace almost central unit vector by a sequence of unit vectors that is ``almost central'' for $A$.

\smnoind
(b) In Proposition \ref{acv}(b)(iii), we only need to check one bimodule (namely, the universal one) in order to verify strong property $(T)$. 

\smnoind
(c) One may wonder if it is possible to check whether a $C^*$-algebra has property $(T)$ by looking at its universal bimodule alone. 
However, this cannot be done using the original formulation of property $(T)$ because there exists a unital $C^*$-algebra $A$ which does not have property $(T)$ but $\FH^A \neq (0)$ (i.e. $A$ has a tracial state). 
Nevertheless, we will show in Theorem \ref{t-T} below that it is possible to do so using an equivalent formulation of property $(T)$.
\end{rem}

\medskip

\section{Kazhdan constants}

\medskip

In this section, we will define and study some Kazhdan constants in the case when $B = A$. Let us start with the following
lemma.

\medskip

\begin{lem}
\label{dec-bimod} Let $H\in \sbm(A)$. If $H_\CC$ is the sub-bimodule
generated by $H^A$ (called the \emph{centrally generated part of
$H$}), then $H_\CC$ is a quasi-standard bimodule (Definition
\ref{def-tq}) and $H_\CC^\bot$ contains no non-zero central vector
for $A$.
\end{lem}
\begin{prf}
Without loss of generality, we may assume that $C := \sph(H^A)$ is
non-empty. 
Let $\mathcal{S} := \{\overline{A\cdot\xi}: \xi\in C\}$ and
$$\mathfrak{M}\ :=\ \{ \mathcal{M}\subseteq \mathcal{S}: K \bot L\ {\rm
for\ any} \ K,L\in \mathcal{M}\}.$$ 
By the Zorn's lemma, there exists
a maximal element $\mathcal {M}_{0}$ in $\mathfrak{M}$ and we put
$H_1 := \bigoplus_{K\in  \mathcal{M}_0} K$. 
Then clearly
$H_1\subseteq H_\CC$ and $H_1^\bot$ contains no non-zero central
vector for $A$. 
Together with the fact $H^A = H_1^A \oplus
(H_1^\bot)^A$, this shows that $H^A = H_1^A \subseteq H_1$ and hence $H_1 =
H_\CC$. 
Finally, for any $\overline{A\cdot\xi}\in \mathcal{M}_0$ with
$\xi\in C$, the functional defined by $\tau(a) := \la a\xi, \xi \ra$
($a\in A$) is a tracial state and $\overline{A\cdot\xi} \cong M_\tau$.
This completes the proof.
\end{prf}

\medskip

Suppose that $H\in\sbm(A)$ and $K$ is a Hilbert subspace of $H$. For
any $Q\in \F(A)$, we set
$$t^A(Q;H,K)\ :=\ \inf\ \left\{ \left(\sum_{x\in
Q} \|x\cdot \xi -\xi\cdot x\|^{2}\right)^{1/2}:  \xi \in \sph(H \ominus K)
\right\}$$ (we use the convention that the infimum over the empty
set is $+\infty$).

\medskip

\begin{lem}
\label{lem-thk} Let $Q\in \F(A)$, $H\in\sbm(A)$ and $K$ be a Hilbert
subspace of $H$. 
Suppose that $H = \bigoplus_{\lambda\in \Lambda}
H_\lambda$ such that $K = \bigoplus_{\lambda\in \Lambda} K_\lambda$
where $K_\lambda := H_\lambda\cap K$.

\smnoind (a) If $\alpha_\lambda$ is a cardinal for any $\lambda\in \Lambda$, and 
if we set $H_0 := \bigoplus_{\lambda\in\Lambda} \left(
\bigoplus_{\alpha_\lambda} H_\lambda \right)$ and $K_0 :=
\bigoplus_{\lambda\in\Lambda} \left(\bigoplus_{\alpha_\lambda}
K_\lambda \right)$, then
\begin{equation}
\label{tq<sum} t^A(Q;H,K)^2\|\zeta\|^2\ \leq\ \sum_{x\in Q}\|x\cdot \zeta
- \zeta \cdot x\|^2 \qquad (\zeta\in H_0 \ominus K_0).
\end{equation}

\smnoind (b) $t^A(Q;H,K) = \inf_{\lambda\in \Lambda}
t^A(Q;H_\lambda,K_\lambda).$

\smnoind (c) $(t^A(Q;H,K))_{Q\in \F(A)}$ is an increasing net and
$\lim_{Q\in \F(A)} t^A(Q;H,K) = 0$ if and only if there exists an
almost central unit vector for $A$ in $H\ominus K$. 
In this case, one can
choose an almost central unit vector $(\xi_i)_{i\in I}$ for $A$ such that
for any $i\in I$, there exists $\lambda_i\in \Lambda$ with $\xi_i\in
H_{\lambda_i} \ominus K_{\lambda_i}$.
\end{lem}
\begin{prf}
(a) For any $\zeta\in \sph(H_0 \ominus K_0)$, we have $\zeta =
(\zeta_{i,\lambda})_{\lambda\in \Lambda; i\in \alpha_\lambda}$
with $\zeta_{i,\lambda}\in H_\lambda \ominus K_\lambda \subseteq
H\ominus K$ and $\sum_{\lambda\in \Lambda} \sum_{i\in
\alpha_\lambda} \|\zeta_{i,\lambda}\|^2 = 1$. Thus,
$$t^A(Q;H,K)^2
\ \leq \ \sum_{\lambda\in \Lambda} \sum_{i\in
\alpha_\lambda}\|\zeta_{i,\lambda}\|^2 \sum_{x\in Q}
\left\|x\cdot \frac{\zeta_{i,\lambda}}{\|\zeta_{i,\lambda}\|} - \frac{\zeta_{i,\lambda}}{\|\zeta_{i,\lambda}\|}\cdot x\right\|^2 \ = \
\sum_{x\in Q} \|x\cdot \zeta - \zeta \cdot x\|^2.$$ 

\smnoind (b) Note that $t^A(Q;H,K) \leq t^A(Q;H_\lambda,K_\lambda)$
for all $\lambda\in \Lambda$ (as $H_\lambda \ominus K_\lambda
\subseteq H\ominus K$). 
For any $\epsilon > 0$, there exists $\xi\in
\sph(H\ominus K)$ such that $\sum_{x\in Q} \|x \cdot \xi - \xi \cdot x\|^2 \leq
t^A(Q;H,K) + \epsilon$. Now, $\xi = (\xi_\lambda)_{\lambda\in
\Lambda}$ with $\xi_\lambda\in H_\lambda \ominus K_\lambda$ and
$\sum_{\lambda\in \Lambda} \|\xi_\lambda\|^2 = 1$. A similar
argument as part (a) implies that there exists $\lambda_0 \in
\Lambda$ such that
$$\sum_{x\in Q} \left\|x \cdot \frac{\xi_{\lambda_0}}
{\|\xi_{\lambda_0}\|} - \frac{\xi_{\lambda_0}}
{\|\xi_{\lambda_0}\|} \cdot x\right\|^2\ \leq\ t^A(Q;H,K) + \epsilon.$$

\smnoind (c) It is clear that $(t^A(Q;H,K))_{Q\in \F(A)}$ is
increasing and that $t^A(Q;H,K) = 0$ for any $Q\in \F(A)$ if there
exists an almost central unit vector for $A$ in $H\ominus K$. Now,
suppose that
$$\sup_{Q\in\F(A)}\inf_{\lambda\in \Lambda} t^A(Q;H_\lambda, K_\lambda)
\ =\ \lim_{Q\in \F(A)} t^A(Q;H,K) \ =\ 0.$$ Then for any $Q\in
\F(A)$ and $\epsilon > 0$, there exists $\lambda_{Q,\epsilon} \in
\Lambda$ and $\xi_{Q,\epsilon} \in \sph(H_{\lambda_{Q,\epsilon}}
\ominus K_{\lambda_{Q,\epsilon}})$ such that $\sum_{x\in Q} \|x\cdot 
\xi_{Q,\epsilon} - \xi_{Q,\epsilon} \cdot x\|^2 < \epsilon^2$. It is easy
to see that $(\xi_{Q,\epsilon})_{(Q,\epsilon)\in \F(A)\times \RP}$
is an almost central unit vector for $A$.
\end{prf}

\medskip

Now, we define three Kazhdan constants: for any $Q\in \F(A)$,
set
$$t^A_u(Q)\ :=\ t^A(Q;\FH, \FH^A), \quad t^A_c(Q)\
:=\ t^A(Q;\FH, \FH_\CC), \quad t^A_s(Q)\ :=\ t^A(Q;\FK, \FK^A)$$
(where $\FH$ and $\FK$ are the universal bimodule and the standard
bimodule respectively) and
$$t^A_u : = \sup_{Q\in \F(A)} t^A_u(Q), \quad t^A_c
:= \sup_{Q\in \F(A)} t^A_c(Q) \quad {\rm as\ well\ as} \quad t^A_s
:= \sup_{Q\in \F(A)} t^A_s(Q).$$

\medskip

\begin{lem}
\label{rem-tq} 
(a) For any $H\in \sbm(A)$, we have $t^A_u(Q) \leq t^A(Q;H,H^A)$ and $t^A_c(Q)\leq t^A(Q;H,H_\CC)$. 
If, in addition, $H$ is quasi-standard, then $t^A_s(Q)\leq t^A(Q;H,H^A)$.

\smnoind
(b) $t^A_u(Q)\leq \min \{t^A_c(Q), t^A_s(Q)\}$.
\end{lem}
\begin{prf}
(a) There are cardinals $\alpha_\tau$ ($\tau\in S(A^{Dou})$) such that $H\cong \bigoplus_{\tau\in S(A^{Dou})}
M_{\alpha_\tau, \tau}$. 
For any $\xi\in \sph(H\ominus H^A)$, we have $\xi = (\xi_\tau)$ where $\xi_\tau\in M_{\alpha_\tau, \tau} \ominus M_{\alpha_\tau, \tau}^A$. 
By Inequality \eqref{tq<sum}, we have
$$t^A_u(Q)^2\|\xi_\tau\|^2\ \leq\ \sum_{x\in Q}\|x\cdot \xi_\tau - \xi_\tau \cdot x\|^2$$
and so
$$t^A_u(Q)^2\ \leq\ \sum_{\tau\in S(A^{Dou})}
\sum_{x\in Q}\|x\cdot \xi_\tau - \xi_\tau \cdot x\|^2\ =\ \sum_{x\in Q}\|x\cdot \xi - \xi \cdot x\|^2$$
(as $\sum_{\tau\in S(A^{Dou})} \|\xi_\tau\|^2 = \|\xi\|^2 = 1$). 
Thus, we have $t^A_u(Q) \leq t^A(Q;H,H^A)$.
The arguments for the other two inequalities are similar. 

\smnoind 
(b) $t^A_u(Q)\leq t^A_c(Q)$ because $\FH^A \subseteq \FH_\CC$ and 
$t^A_u(Q)\leq t^A_s(Q)$ because of part (a). 
\end{prf}

\medskip

\begin{thm}
\label{t-T} (a) The following statements are equivalent.
\begin{enumerate}
\item[\rm i.] $t^A_u > 0$.
\item[\rm ii.] $A$ has strong property $(T)$.
\item[\rm iii.] There exists $(Q, \delta) \in \F(A)\times \RP$
such that for any $\xi\in V_\FH(Q,\delta)$, we have $P^A_\FH(\xi)\neq 0$.
\end{enumerate}

\smnoind (b) The following statements are equivalent.
\begin{enumerate}
\item[\rm i.] $t^A_s > 0$.
\item[\rm ii.] For any $\epsilon > 0$, there exists
$(Q, \delta) \in \F(A)\times \RP$ such that for any quasi-standard
bimodule $H$ and any $\xi\in V_H(Q,\delta)$,
we have $\|\xi - P^A_H(\xi)\| < \epsilon$.
\item[\rm iii.] There exists $(Q, \delta) \in \F(A)\times \RP$
such that for any $\xi\in V_{\FK}(Q,\delta)$,
we have $P^A_{\FK}(\xi)\neq 0$.
\end{enumerate}

\smnoind (c) The following statements are equivalent.
\begin{enumerate}
\item[\rm i.] $t^A_c > 0$.
\item[\rm ii.] $A$ has property $(T)$.
\item[\rm iii.] There is $(Q,\delta)\in \F(A)\times \RP$
such that $V_H(Q,\delta)\cap H_\CC^\bot = \emptyset$
for any $H\in \sbm(A)$.
\item[\rm iv.] There exists $(Q,\delta)\in \F(A)\times \RP$
such that $V_\FH(Q,\delta)\cap \FH_\CC^\bot = \emptyset$.
\end{enumerate}

\smnoind (d) $t^A_u > 0$ if and only if both $t^A_c > 0$ and $t^A_s
> 0$.
\end{thm}
\begin{prf}
(a) (i)$\Rightarrow$(ii). There exists $Q\in \F(A)$ with $t^A_u(Q) >
0$. Let $m$ be the number of elements in $Q$ and $\delta =
\frac{t^A_u(Q)\epsilon}{\sqrt{m}}$. For any $H\in \sbm(A)$ and
$\tau\in S(A^{Dou})$, there is a cardinal $\alpha_\tau$ such that $H
= \bigoplus _{\tau \in S(A^{Dou})} M_{\alpha_\tau, \tau}$
($\alpha_\tau$ can be zero). Pick any $\xi\in V_H(Q,\delta)$ and
consider $\xi'' = \xi - P_H^A(\xi)\in (H^A)^\bot$. Since $\xi'' =
(\zeta_\tau)_{\tau\in S(A)}$ where $\zeta_\tau \in (M_{\alpha_\tau,
\tau}^A)^\bot$, we have, by Inequality
\eqref{tq<sum},
$$\|\xi^{''}\|^{2}
\ \leq\ t^A_u(Q)^{-2}\sum_{x\in Q} \| x\cdot \xi^{''} -\xi^{''}\cdot x\|^{2} \ =
\ t^A_u(Q)^{-2}\underset{x\in Q}{\sum} \| x\cdot \xi -\xi \cdot x\|^{2} \ <\
\epsilon^2.$$

\smnoind (ii)$\Rightarrow$(iii). By taking $\epsilon = 1/2$, we see
that Statement (iii) holds.

\smnoind (iii)$\Rightarrow$(i). Suppose on the contrary that $t^A_u(Q) = 0$. Then
there exists $\xi \in \sph((\FH^A)^\bot)$ with $\sum_{x\in Q} \|
x\cdot \xi -\xi \cdot x\|^{2} < \delta^2$. Hence, $\xi\in V_{\FH}(Q,\delta)$ and
so $P_{\FH}^A(\xi)\neq 0$ which contradicts the fact that $\xi\in
(\FH^A)^\bot$. 

\smnoind (b) The proof of this part is essentially the same as that
of part (a) with $\FH$ and $t^A_u$ being replaced by $\FK$ and
$t^A_s$ respectively.

\smnoind (c) (i)$\Rightarrow$(ii). Let $Q\in \F(A)$ such that
$t^A_c(Q) > 0$. Suppose that $A$ does not have property $(T)$. There
exists $H\in \sbm(A)$ that contains an almost central unit vector
$(\xi_i)$ for $A$ but $H^A = \{0\}$. 
Hence, $H_\CC^\bot = H$, and 
\begin{equation*}
\label{inf=0} 
t^A(Q;H, H_\CC)\ =\ \inf\ \left\{ \left(\sum_{x\in Q} \| x\cdot \xi -\xi
\cdot x\|^{2}\right)^{1/2}:  \xi \in \sph(H) \right\}\ =\ 0.
\end{equation*}
Now Lemma \ref{rem-tq}(a) gives the contradiction that $t^A_c(Q)  = 0$.

\smnoind (ii)$\Rightarrow$(i). Suppose on the contrary that $t^A_c(Q) = 0$ for any
$Q\in \F(A)$. There exists, by Lemma \ref{lem-thk}(c), an almost
central unit vector for $A$ in $\FH_\CC^\bot$ which contradicts the fact that $A$
has property $(T)$ (because of Lemma \ref{dec-bimod}).

\smnoind (i)$\Rightarrow$(iii). Let $Q\in \F(A)$ such that $t^A_c(Q)
> 0$ and $0 < \delta < \frac{t^A_c(Q)}{\sqrt{m}}$ where $m$ is the
number of elements in $Q$. 
By Lemma \ref{rem-tq}(a), we have
$t^A_c(Q) \leq t^A_c(Q;H, H_\CC)$. 
Thus, 
$$\sum_{x\in Q}\|x\cdot \zeta - \zeta \cdot x\|^2\ \geq\ t_c^A(Q)\ >\ m\delta^2$$ for any $\zeta\in \sph(H_\CC^\bot)$.
Suppose that there exists $\xi \in V_H(Q,\delta)\cap H_\CC^\bot$.
Then $m\delta^2 < \sum_{x\in Q}\|x\cdot \xi - \xi \cdot x\|^2 < m\delta^2$
which is absurd.

\smnoind (iii)$\Rightarrow$(iv). This is obvious.

\smnoind (iv)$\Rightarrow$(i). Suppose on the contrary that $t^A_c(Q) = 0$. 
Then
there exists $\xi \in \sph(\FH_\CC^\bot)$ with $\sum_{x\in Q} \|
x\cdot \xi -\xi \cdot x\|^{2} < \delta^2$ which gives the contradiction that
$\xi\in V_{\FH}(Q,\delta)\cap \FH_\CC^\bot$. 

\smnoind (d) If $t^A_u > 0$, then $t^A_s > 0$ and $t^A_c > 0$ (by Lemma \ref{rem-tq}(b)). 
Conversely,
suppose that $t^A_u = 0$. Then by Lemma \ref{lem-thk}(c), there
exists an almost central unit vector $(\xi_i)_{i\in I}$ for $A$ in 
$$\FH \ominus \FH^A\ =\ (\FH \ominus \FH_\CC) \oplus (\FH_\CC \ominus
\FH^A).$$ 
Let $\eta_i \in \FH \ominus \FH_\CC$ and $\zeta_i \in
\FH_\CC \ominus \FH^A$ be the corresponding components of $\xi_i$.
Then either $\eta_i \nrightarrow 0$ or $\zeta_i \nrightarrow 0$.
Therefore, by rescaling, there exists an almost central unit vector for
$A$ in either $\FH \ominus \FH_\CC$ or $\FH_\CC \ominus \FH^A =
\FH_\CC \ominus \FH_\CC^A$. In the first case, we have $t^A_c = 0$
(by Lemma \ref{lem-thk}(c)). In the second case, we have $t^A_s \leq
\sup_{Q\in \F(A)} t^A(Q; \FH_\CC, \FH_\CC^A) = 0$ (by Lemma
\ref{rem-tq}(a), Lemma \ref{dec-bimod} and Lemma \ref{lem-thk}(c)).
\end{prf}

\medskip

Part (a) of the above theorem tells us that in order to show that
$A$ has strong property $(T)$, it suffices to verify a weaker
condition than that of Definition \ref{def-sT} for just the
universal bimodule $\FH$.

\medskip

\begin{rem}
\label{rem-equ-t} (a) The argument of Theorem \ref{t-T}(a), together
with Lemma \ref{rem-tq}(a),  tell us that for any $Q\in \F(A)$,
$\delta > 0$ and $H\in \sbm(A)$, if $\xi\in V_H(Q,\delta)$, then
$$t^A_u(Q) \|\xi - P_H^A(\xi)\|\ <\ \delta \sqrt{m}$$ 
(where $m$ is the number of elements in $Q$).

\smnoind (b) The argument of Theorem \ref{t-T}(c), together with
Lemma \ref{rem-tq}(a), tell us that if $t^A_c(Q) > 0$, for any $Q\in \F(A)$, $H\in \sbm(A)$ and $\delta \in (0, \frac{t^A_c(Q)}{\sqrt{m}})$, we have $V_H(Q,\delta)\cap H_\CC^\bot =
\emptyset$.

\smnoind
(c) The gap between property $(T)$ and strong property $(T)$ is represented by the gap between $t_c^A$ and $t_u^A$ or equivalently between $\FH^A$ and $\FH_\CC$. 
Note that in the case of a locally compact group $G$, such a gap does not exist because the set of $G$-invariant vectors defines a subrepresentation. 
\end{rem}

\medskip

By Theorem \ref{t-T} and Lemma \ref{lem-thk}(c), we have the following corollary. 

\medskip

\begin{cor}
\label{cor-main-t} (a) $A$ has property $(T)$ (respectively, strong
property $(T)$) if and only if there is no almost central unit
vector for $A$ in $\FH_\CC^\bot$ (respectively, in $(\FH^A)^\bot$).

\smnoind (b) $A$ has strong property $(T)$ if and only if $A$ has
property $(T)$ and $t^A_s > 0$.
\end{cor}

\medskip

Note that one can also obtain part (b) of the above corollary by
using a similar argument as that of \cite[Proposition 1]{CJ}.

\medskip

\section{Some permanence properties}

\medskip

In this section, we study the permanence properties for property
$(T)$ and strong property $(T)$. First of all, we have the following
lemma which implies that the quotient of any pair having property
$(T)$ (respectively, strong property $(T)$) will have the same
property. Since the proof is direct, we will omit it.

\medskip

\begin{lem}
\label{quot-t} Let $A_1$ and $A_2$ be two unital $C^*$-algebras and let $B_1\subseteq A_1$ and $B_2\subseteq A_1$ be $C^*$-subalgebras containing the identities of $A_1$ and $A_2$ respectively. 
Suppose that $\varphi: A_1\rightarrow A_2$ is a unital $*$-homomorphism such that $B_2\subseteq \varphi(B_1)$.  
If $(A_1, B_1)$ has property $(T)$ (respectively, strong property $(T)$), then so does
$(A_2, B_2)$.
\end{lem}

\medskip

\begin{lem}
\label{sum} Let $A_1$, $A_2$, $B_1$ and $B_2$ be the same as Lemma \ref{quot-t}. 
If both $(A_1, B_1)$ and $(A_2, B_2)$ have property
$(T)$ (respectively, strong property $(T)$), then so does $(A_1\oplus
A_2, B_1\oplus B_2)$.
\end{lem}
\begin{prf}
The statement for property $(T)$ is well known and we will only show
the case for strong property $(T)$. Suppose that $H\in
\sbm(A_1\oplus A_2)$ and $e = (1_{A_1}, 0) \in A_1\oplus A_2$. Then
$H = \bigoplus_{k,l =1}^2 H_{kl}$ where $H_{kl}$ is a non-degenerate
Hilbert $A_k$-$A_l$-bimodule. Suppose that $(\xi_i)_{i\in I}$ is an
almost central unit vector in $H$ for $A_1\oplus A_2$ and $\xi_i =
\sum_{k,l =1}^2 \xi^{kl}_i$ where $\xi^{kl}_i\in H_{kl}$. Then
$$\left\|\xi^{12}_i\right\|^2 + \left\|\xi^{21}_i\right\|^2\
=\ \left\|e \cdot \xi_i - \xi_i \cdot e\right\|^2\ \rightarrow\ 0.$$
If $\left\|\xi^{22}_i\right\| \rightarrow 0$, then we can assume that $\|\xi^{11}_i\| > 1/2$ ($i\in I$), and $\left(\frac{\xi^{11}_i}{\|\xi^{11}_i\|}\right)_{i\in I}$ is an
almost central unit vector for $A_1$ in $H_{11}$. 
In this case, for any $\epsilon >
0$, there is $i_0\in I$ such that $\left\|\xi^{11}_i -
P_{H_{11}}^{B_1}\left(\xi^{11}_i\right)\right\| <
\frac{\epsilon}{2}$ ($i\geq i_0$), which implies that
$$\left\|\xi_j - P_{H}^{B_1\oplus B_2}(\xi_j)\right\|\
\leq\ \sqrt{\left\|\xi^{11}_i - P_{H_{11}}^{B_1}
(\xi^{11}_i)\right\|^2 + \left\|\xi_i^{12}\right\|^2
+ \left\|\xi_i^{21}\right\|^2 + \left\|\xi_i^{22}\right\|^2}\ <\ \epsilon$$
when $j$ is large enough. The same conclusion holds if
$\left\|\xi_i^{11}\right\| \rightarrow 0$. 
We consider now the case
when $\left\|\xi_i^{11}\right\| \nrightarrow 0$ and
$\left\|\xi^{22}_i\right\| \nrightarrow 0$. 
There exist a
constant $\kappa > 0$ as well as subnets
$\left(\xi^{11}_{i_k}\right)_{k\in J_1}$ and
$\left(\xi^{22}_{i_l}\right)_{l\in J_2}$ such that
$\left\|\xi^{11}_{i_k}\right\|, \left\|\xi^{22}_{i_l}\right\| \geq
\kappa$ for every $k\in J_1$ and $l\in J_2$. One can show easily
that $\left(\frac{\xi^{11}_{i_k}}{\|\xi^{11}_{i_k}\|}\right)_{k\in
J_1}$ and
$\left(\frac{\xi^{22}_{i_l}}{\|\xi^{22}_{i_l}\|}\right)_{l\in J_2}$
are almost central unit vectors for $A_1$ and $A_2$ in $H_{11}$ and
$H_{22}$ respectively. 
Thus, for any $\epsilon > 0$, one can find $i_0\in I$ such that
$$\left\|\xi^{11}_{i_0} - P_{H_{11}}^{B_1}\left(\xi^{11}_{i_0}\right)\right\|
\ <\ \frac{\left\|\xi^{11}_{i_0}\right\|\epsilon}{2}, \qquad
\left\|\xi^{22}_{i_0} -
P_{H_{22}}^{B_2}\left(\xi^{22}_{i_0}\right)\right\| \ <\
\frac{\left\|\xi^{22}_{i_0}\right\|\epsilon}{2}$$ and
$\left\|\xi^{12}_{i_0}\right\|^2 + \left\|\xi^{21}_{i_0}\right\|^2 \
<\ \frac{\epsilon^2}{2}$. Consequently,
$$\|\xi_{i_0} - P_{H}^{B_1\oplus B_2}(\xi_{i_0})\|\
\ \leq\\ \sqrt{\|\xi^{11}_{i_0} -
P_{H_{11}}^{B_1}(\xi^{11}_{i_0})\|^2 + \|\xi^{22}_{i_0} -
P_{H_{11}}^{B_1}(\xi^{11}_{i_0})\|^2 + \|\xi_{i_0}^{12}\|^2 +
\|\xi_{i_0}^{21}\|^2} \ <\ \epsilon.$$ In any case, $(A_1\oplus A_2,
B_1\oplus B_2)$ has strong property $(T)$ because of Proposition 
\ref{acv}(b)(iii).
\end{prf}

\medskip

Our next task is to consider tensor products and crossed products. 
Let
us first recall the following useful terminology of co-rigidity from
\cite[Remark 19]{Bek-T}. We will also introduce a stronger version
of co-rigidity corresponding to strong property $(T)$.

\medskip

\begin{defn}
\label{def-corig} 
The pair $(A,B)$ is said to be

\smnoind (a) \emph{co-rigid} if there exists $(Q,\beta)\in
\F(A)\times \RP$ such that for any $H\in \sbm(A)$ with $V_H(Q,\beta)
\cap H^B \neq \emptyset$, we have $H^A\neq \{0\}$.

\smnoind (b) \emph{strongly co-rigid} if for any $\gamma > 0$, there
exists $(Q,\delta)\in \F(A)\times \RP$ such that for any $H\in
\sbm(A)$ and any $\xi\in V_H(Q,\delta) \cap H^B$, we have
$\left\|\xi - P^A_H(\xi)\right\| < \gamma$.
\end{defn}

\medskip

The idea of the following result comes from \cite[2.3]{Ana}.

\medskip

\begin{prop}
\label{T-corig} Suppose that $B$ has strong property $(T)$.

\smnoind (a) $A$ has property $(T)$ if and only if $(A,B)$ is
co-rigid.

\smnoind (b) $A$ has strong property $(T)$ if and only if $(A,B)$ is
strongly co-rigid.
\end{prop}
\begin{prf}
(a) The sufficiency is clear and we will only show the necessity.
Let $(Q,r)\in \F(A)\times \RP$ be the pair satisfying the condition
in Definition \ref{def-corig}(a). 
Suppose that $(F,s)\in \F(B)\times
\RP$ is the strong Kazhdan's pair for $(B, B, \alpha)$ where $\alpha = \min \{\frac{r}{8M}, \frac{1}{2}\}$ and $M = \max \{\|a\|: a\in Q\}$. 
Put $E = Q\cup F$ and $t =
\min \{ \frac{r}{4}, s\}$. Assume that $H\in \sbm(A)$ with $\xi\in
V_H(E,t)$. As $\xi\in V_H(F,s)$, one has $\left\| \xi -
P^B_H(\xi)\right\| < \alpha$ and $\left\|P^B_H(\xi)\right\| \geq
\frac{1}{2}$. If $\eta =
\frac{P^B_H(\xi)}{\left\|P^B_H(\xi)\right\|}$, then we have, for any
$a\in Q$,
$$\|a\cdot \eta - \eta\cdot a\|
\ \leq \ \frac{\|a\cdot \xi - \xi\cdot a\| + 2\|a\|\left\| \xi -
P^B_H(\xi)\right\|}{\left\|P^B_H(\xi)\right\|} \ < \ 2t +
\frac{r\|a\|}{2M} \ \leq \ r.$$ Thus, $\eta\in V_H(Q,r)\cap H^B$ and
$H^A\neq \{0\}$.

\smnoind (b) Again, we only need to show the necessity. For any
$\epsilon > 0$, let $(Q,r)\in \F(A)\times \RP$ be the pair
satisfying Definition \ref{def-corig}(b) for $\gamma =
\frac{\epsilon}{2}$. 
Take a strong Kazhdan's pair $(F,s)\in \F(B)\times \RP$ for $(B,B,\alpha)$ where $\alpha = \min
\{\frac{r}{8M}, \frac{1}{2}, \frac{\epsilon}{2}\}$ and $M = \max \{\|a\|: a\in Q\}$. 
If $E$ and $t$
are as in the argument of part (a), then for any $\xi\in V_H(E,t)$,
we have $\eta = \frac{P^B_H(\xi)}{\left\|P^B_H(\xi)\right\|}\in
V_H(Q,r)\cap H^B$ which implies that
$$\left\| P^B_H(\xi) - P^A_H(\xi)\right\|
\ < \ \frac{\left\| P^B_H(\xi)\right\|\epsilon}{2} \ \leq \
\frac{\epsilon}{2}$$ (note that $H^A \subseteq H^B$). Since $\left\|
\xi - P^B_H(\xi)\right\| < \frac{\epsilon}{2}$ as well, we see that
$(E,t)\in \F(A)\times \RP$ is a strong Kazhdan's pair for $(A, B, \epsilon)$.
\end{prf}

\medskip

We do not know whether $B$ having property $(T)$ and $(A,B)$ being
co-rigid will imply that $A$ has property $(T)$. If it is the case,
then the statement in Theorem \ref{ten}(a) below concerning property $(T)$
can be improved and a similar statement as Theorem \ref{st=t+q} below for
property $(T)$ will also hold.

\medskip

The first application of the above proposition is the following theorem. 
Notice that unlike the case of type $\rm II_1$-factors (see \cite[2.5]{Ana}), the fact that $B\otimes_{\rm max} D$ having property
$(T)$ (or strong property $(T)$) will not imply both $B$ and $D$
to have property $(T)$ (respectively, strong property $(T)$), but at least one of them have property $(T)$ (respectively, strong property $(T)$).

\medskip

\begin{thm}
\label{ten} Let $B$ and $D$ be two unital $C^*$-algebras, $A =
B\otimes_{\rm max} D$ and $A_0 = B\otimes_{\rm min} D$.

\smnoind (a) If $B$ has strong property $(T)$ and $D$ has property
$(T)$ (respectively, strong property $(T)$), then $A$ has property
$(T)$ (respectively, strong property $(T)$).

\smnoind (b) If there is no almost central unit vector for $D$ in any $K\in \sbm(D)$, then $A$ has strong property $(T)$.

\smnoind (c) Suppose that there exists an almost central unit vector $(\eta_j)_{j\in J}$ for
$D$ in some $K\in \sbm(D)$. 
If $A_0$ has property $(T)$ (respectively,
strong property $(T)$), then so does $B$.

\smnoind (d) If $A_0$ has property $(T)$ (respectively, strong
property $(T)$), then either $B$ or $D$ has property $(T)$
(respectively, strong property $(T)$).
\end{thm}
\begin{prf}
(a) We show the statement for strong property $(T)$ first. Suppose
that $\alpha > 0$ and $(F,r)\in \F(D)\times \RP$ is the strong Kazhdan's pair for $(D,D,\alpha)$. 
Let $Q := 1\otimes F \in \F(A)$ and $H\in \sbm(A)$. If $H^B \neq
\{0\}$, then $H^B\in \sbm(D)$ under the canonical multiplications.
For any $\xi\in V_H(Q, r) \cap H^B = V_{H^B}(F,r)$, we have
$$\|\xi - P^A_{H}(\xi)\|\ \leq\ \|\xi - P^D_{H^B}(\xi)\|\ <\ \gamma$$
(note that $(H^B)^D =H^A$). Thus, $(A,B)$ is strongly co-rigid and
we can apply Proposition \ref{T-corig}(b). The proof for the case of
property $(T)$ is similar.

\smnoind (b) In this case, there is no almost central
unit vector for $A$ in any $H\in \sbm(A)$ and we can apply Proposition \ref{acv}.

\smnoind (c) We will establish the statement for strong property $(T)$ and the
statement for property $(T)$ is similar (and easier). 
Suppose $H\in \sbm(B)$ and there exists an almost central
unit vector $(\xi_i)_{i\in I}$ for $B$ in $H$. 
Then $(\xi_i\otimes
\eta_j)_{(i,j)\in I\times J}$ is an almost central unit vector for $A_0$
in $H\otimes K$. For any $\epsilon > 0$, there exists, by
Proposition \ref{acv}(b)(ii),  $(i_0,j_0)\in I\times J$ such that
$$\left\|\xi_{i_0}\otimes \eta_{j_0} - P^{A_0}_{H\otimes K}
(\xi_{i_0}\otimes \eta_{j_0})\right\|\ <\ \epsilon.$$
Let $\zeta = P^{A_0}_{H\otimes K}(\xi_{i_0}\otimes \eta_{j_0})\in
(H\otimes K)^{A_0}$ and $\varphi\in K^*$ be defined by
$\varphi(\eta') = \la \eta', \eta_{j_0}\ra$. For any $b\in B$ and
any $\xi\in H$, we have
$$\la b\cdot (\id\otimes \varphi)(\zeta), \xi \ra\
=\ \la (b\otimes 1)\cdot\zeta, \xi\otimes \eta_{j_0} \ra\
=\ \la \zeta\cdot (b\otimes 1), \xi\otimes \eta_{j_0} \ra\
=\ \la (\id\otimes \varphi)(\zeta)\cdot b, \xi \ra.$$
Therefore, $(\id\otimes \varphi)(\zeta)\in H^B$ and
$$\|\xi_{i_0} - (\id\otimes \varphi)(\zeta)\|\
=\ \left\|(\id\otimes \varphi)\left(\xi_{i_0}
\otimes \eta_{j_0} - P^A_{H\otimes K}(\xi_{i_0}
\otimes \eta_{j_0})\right)\right\|\ <\ \epsilon.$$
This shows that $B$ has strong property $(T)$ (by Proposition \ref{acv}(b)(iii)).

\smnoind (d) 
If there is no almost central unit vector for $D$ in any $K\in \sbm(D)$, then $D$ has strong property $(T)$ (by definition). 
Otherwise, we can apply part (b).
\end{prf}

\medskip

Next, we will consider crossed product of $C^*$-algebras by actions of discrete groups. 
The idea of which comes from \cite[4.6]{Ana}. 
Again, unlike the case of type $\rm II_1$-factors, even if $B\times_\alpha \Gamma$ has strong property $(T)$, this will not imply that $\Gamma$ has
property $(T)$ (notice that if $\alpha$ is trivial and any element in $\sbm(B)$ does not contain an almost central unit vector, then $B\times_\alpha \Gamma$ will have strong property $(T)$ whether or not $\Gamma$ have property $(T)$).

\medskip

\begin{thm}
\label{st=t+q} Let $B$ be a unital $C^*$-algebra with an action
$\alpha$ by a discrete group $\Gamma$ and $A = B\times_\alpha
\Gamma$. If $\Gamma$ has property $(T)$, then $(A,B)$ is strongly
co-rigid. Consequently if $B$ has strong property $(T)$ and $\Gamma$
has property $(T)$, then $A$ has strong property $(T)$ (and so does
$B\times_{\alpha, r} \Gamma$).
\end{thm}
\begin{prf}
As $\Gamma$ has property $(T)$, for any $\epsilon > 0$, there exists
$(F,\delta)\in \F(\Gamma)\times \RP$ such that for any $(K,\pi)\in
{\rm Rep}(\Gamma)$ and any $\eta\in V_\pi(F,\delta)$, one has
$\left\|\eta - P^{\Gamma}_K(\eta)\right\| < \epsilon$ (by \cite[Theorem 1.2(b2)]{Jol}). 
Let $\mu:
B\rightarrow A$ and $u:\Gamma \rightarrow A$ be the canonical maps.
For any $H\in \sbm(A)$, we define a representation $\pi: \Gamma
\rightarrow \CL(H^B)$ by $\pi(t)\xi = u_t\cdot \xi \cdot u_t^*$ ($t\in \Gamma, \xi\in H^B$) which is well defined because for any $b\in B$, 
$$\mu(b)\cdot \pi(t)\xi
\ =\ u_t\mu(\alpha_{t^{-1}}(b))\cdot \xi \cdot u_t^*
\ =\ u_t \cdot \xi \cdot \mu(\alpha_{t^{-1}}(b))u_t^*
\ =\ \pi(t)\xi\cdot \mu(b).$$ 
Moreover, it is easy to check that $(H^B)^{\Gamma} = H^A$. Thus,
if $\eta\in V_H(u(F), \delta)\cap H^B = V_\pi(F, \delta)$, then
$\left\|\eta - P^A_H(\eta)\right\| = \left\|\eta -
P^{\Gamma}_{H^B}(\eta)\right\| < \epsilon$. This shows that $(A,B)$
is strongly co-rigid. The last statement follows from Proposition
\ref{T-corig}(b).
\end{prf}

\medskip

\section{Some examples of strong property $(T)$}

\medskip

Our first example is finite dimensional $C^*$-algebras. It is easy
to see that any element in $\sbm(M_n(\mathbb{C}))$ has a non-zero
central vector and so $M_n(\mathbb{C})$ has property $(T)$. In fact,
$M_n(\mathbb{C})$ also has strong property $(T)$ but a bit more
argument is needed to establish this fact.

\medskip

\begin{eg}
\label{matrix} $A = M_n(\mathbb{C})$ has a unique tracial state $\tau$ and hence $\FK = M_\tau$. 
If $t^A_s = 0$, there exists an almost central unit vector in $(M_\tau^A)^\bot$ for $A$ (by Lemma \ref{lem-thk}(c)). 
As $\sph((M_\tau^A)^\bot)$ is compact, this implies the existence of a central unit vector $\eta$ for $A$ in $(M_\tau^A)^\bot$ but such $\eta \in M_n(\mathbb{C})$ should be an element of $\mathbb{C} 1 = M_\tau^A$ which is absurd. 
Now by Corollary \ref{cor-main-t}(b) and Lemma \ref{sum}, we see that any finite dimensional
$C^*$-algebra has strong property $(T)$.
\end{eg}

\medskip

%It is also possible to show directly that $M_n(\mathbb{C})$ has
%strong property $(T)$ from the definition (by using the fact that
%any representation of $M_n(\mathbb{C})^{Dou}$ is of the form
%$HS(\mathbb{C}^n)\otimes K$ for some Hilbert space $K$). However, it
%is comparatively harder to do it this way.
By the argument of \cite[Remark 17]{Bek-T}, we know that if $A$ does not have tracial
states, then there is no almost central unit vector for $A$ in any $H\in \sbm(A)$. This,
together with Proposition \ref{acv}(b), give the following result (cf. \cite[Remark
17]{Bek-T}).

\begin{prop}
\label{no-tr-T} If $A$ has no tracial states, then $A$ has strong property $(T)$.
\end{prop}

\medskip

If $H$ is any Hilbert space, $\CL(H)$ has strong property $(T)$ (because of Example \ref{matrix} and
Proposition \ref{no-tr-T}) and so if $B\subseteq \CL(H)$ is any unital $C^*$-subalgebra, then
$(\CL(H), B)$ has strong property $(T)$. 
On the other hand, Proposition \ref{no-tr-T}, together with
Example \ref{matrix}, Lemma \ref{sum} and \cite[5.1]{Bro}, imply the following result.

\medskip

\begin{prop}
Suppose that $A$ is separable and amenable. Then $A$ has strong property $(T)$ if and only
if $A = B\oplus C$ where $B$ is finite dimensional and $C$ has no tracial states. 
\end{prop}

\medskip

We end this section with the following analogue of \cite[Theorem 7]{Bek-T}. 

\medskip

\begin{prop}
\label{eq-disc-gp} Let $\Gamma$ be a countable discrete group and
$\Lambda$ be a subgroup of $\Gamma$. The following statements are
equivalent.

\smnoind (i). $(\Gamma, \Lambda)$ has property $(T)$.

\smnoind (ii). $(C^*(\Gamma), C^*(\Lambda))$ has strong property
$(T)$.

\smnoind (iii). $(C^*(\Gamma), C^*(\Lambda))$ has property $(T)$.

\smnoind (iv). $(C_r^*(\Gamma), C_r^*(\Lambda))$ has strong property
$(T)$.

\smnoind (v). $(C_r^*(\Gamma), C_r^*(\Lambda))$ has property $(T)$.
\end{prop}
\begin{prf}
It is clear that (ii)$\Rightarrow$(iii)$\Rightarrow$(v) and
(ii)$\Rightarrow$(iv)$\Rightarrow$(v) (by Lemma \ref{quot-t}). The
implication (v)$\Rightarrow$(i) was proved in \cite{Bek-T}. It
remains to show that (i)$\Rightarrow$(ii). As $(\Gamma, \Lambda)$
has property $(T)$, for any $\alpha > 0$, there exists $(Q,\beta)\in
\F(\Gamma)\times \RP$ such that for any unitary representation $\pi:
\Gamma \rightarrow \CL(K)$ and any $\xi\in V_\pi(Q,\beta)$, one has
$\|\xi - P^\Lambda_K(\xi)\| < \alpha$ (by \cite[Theorem 1.2(b2)]{Jol}).
Consider $\Gamma \subseteq C^*(\Gamma)$. For any $H\in
\sbm(C^*(\Gamma))$, one can define a unitary representation $\pi_H:
\Gamma \rightarrow \CL(H)$ by $\pi_H(t) \eta = t\cdot \eta \cdot
t^{-1}$ ($\eta\in H$). If $\xi \in V_H(Q,\beta)\subseteq
V_\pi(Q,\beta)$, we have $\left\|\xi -
P^{C^*(\Lambda)}_H(\xi)\right\| = \left\|\xi -
P^\Lambda_H(\xi)\right\| < \alpha$ as required.
\end{prf}

\medskip

\section{Property $(T)$ for the unitary group of a $C^*$-algebra}

\medskip

It was shown in \cite{Pat} that a unital $C^*$-algebra is amenable
if and only if its unitary group under the weak topology is
amenable. 
Motivated by this result as well as by Lemma \ref{F(U)}, we
study in this section, the relation between property $(T)$ and strong property $(T)$ of a unital $C^*$-algebra $A$ and certain $(T)$-type properties of
its unitary group.

\medskip

\begin{rem}
\label{rep>ur} We denote by $\Phi_A: U(A) \rightarrow U(A^{Dou})$
the group homomorphism $u \mapsto u\otimes (u^*)^{\rm op}$. 
Note that any
$K\in \sbm(A)$ defines a non-degenerate $*$-representation $\mu_K$
of $A^{Dou}$ and $\pi_K = \mu_K\circ \Phi_A$ is a unitary
representation of $U(A)$.
\end{rem}

\medskip

The proof of the following result is more or less the same as the
argument for the equivalence of (b2) and (b3) in \cite{Jol}.

\medskip

\begin{prop}
(a) $(A,B)$ has property $(T)$ if and only if for any $K\in
\sbm(A)$, the weak containment of the trivial representation
$1_{U(A)}$ of $U(A)$ in the unitary representation $\pi_K$ will
imply that $\pi_K\!\mid\!_{U(B)}$ contains the trivial
representation $1_{U(B)}$ of $U(B)$.

\smnoind (b) $(A,B)$ has strong property $(T)$ if and only if for
any net $(f_i)_{i\in I}$ in $S(A^{Dou})$ with $(f_i\circ
\Phi_A)_{i\in I}$ converges pointwisely to $1_{U(A)}$, one has
$(f_i\circ\Phi_A\!\mid\!_{U(B)})_{i\in I}$ converges uniformly to
$1_{U(B)}$ on $U(B)$.
\end{prop}
\begin{prf}
(a) Note that for any $K\in \sbm(A)$ and any $\xi\in \sph(K)$, we
have
$$\|t \cdot \xi -\xi \cdot t\|^2 \ = \ 2 - 2 {\rm Re}\la \pi_K(t)\xi,
\xi\ra \qquad (t\in U(A)).$$
Now, this follows almost directly from Proposition \ref{acv}(a).

\smnoind (b) $\Rightarrow)$. 
If $(H_i, \mu_i, \xi_i)$ is the GNS
representation of $f_i$, then $H_i\in \sbm(A)$. 
For any $\epsilon >
0$, let $(Q,\beta)\in \F(U(A))\times \RP$ be a strong Kazhdan's pair for $(A,B, \epsilon/2)$ (see Lemma \ref{F(U)}). 
By the assumption of $(f_i)$, there
exists $i_0\in I$ such that for any $i \geq i_0$, we have
$\sup_{u\in Q} \abs{f_i(\Phi_A(u)) - 1} < \beta^2/2$. Thus,
$$\|u\cdot \xi_i - \xi_i \cdot u\|
\ =\ \|\mu_i(u\otimes (u^*)^{\rm op}) \xi_i - \xi_i\| \ =\
\sqrt{2{\rm Re}(1- f_i(\Phi_A(u)))} \ <\ \beta$$ for any $u\in Q$
and so $\xi_i \in V_{H_i}(Q,\beta)$. Therefore, $\|\xi_i -
P^B_{H_i}(\xi_i)\| < \epsilon/2$ and for any $v\in U(B)$,
$$
\abs{f_i(\Phi_A(v)) -1} \ = \ \abs{\la v\cdot \xi_i \cdot v^* -
\xi_i, \xi_i \ra} \ \leq \ \left\|v\cdot \xi_i \cdot v^* - v \cdot
P^B_{H_i}(\xi_i)\cdot v^*\right\| + \left\|P^B_{H_i}(\xi_i) -
\xi_i\right\| \ < \ \epsilon.
$$

\noindent $\Leftarrow)$. Suppose on the contrary that $(A,B)$ does
not have strong property $(T)$. Let $I := \F(U(A)) \times \RP$. There
exists $\alpha_0 > 0$ such that for any $i = (F,\epsilon)\in I$, one
can find $H_i\in \sbm(A)$ and $\xi_i \in V_{H_i}(F, \epsilon)$ with
$\left\|\xi_i - P^B_{H_i}(\xi_i)\right\| > \alpha_0$. For every such
$H_i$, let $\pi_{H_i}$ be as in Remark \ref{rep>ur}. By
\cite[2.2]{Jol}, we see that there exists $v_i\in U(B)$ such that
\begin{equation}
\label{jol2.2} \|\pi_{H_i}(v_i)\xi_i - \xi_i\|\ >\ \alpha_0
\end{equation}
(note that $H_i^B = H_i^{U(B)}$). On the other hand, for any $i\in
I$, we define $f_i\in S(A^{Dou})$ by $$f_i(a\otimes b^{\rm op}) =
\la a\cdot \xi_i \cdot b, \xi_i\ra \qquad (a,b\in A).$$ For any $i =
(Q,\beta)\in I$, we have $\xi_i \in V_{H_i}(Q, \beta)$ and thus,
$$\sup_{u\in Q}\|\pi_{H_i}(u)\xi_i - \xi_i\|\ <\ \beta.$$ This shows
that $f_i\circ \Phi_A$ converges pointwisely to $1_{U(A)}$.
Therefore, by the hypothesis, \mbox{$f_i\circ
\Phi_A\!\mid\!_{U(B)}$} converges uniformly to $1_{U(B)}$ on $U(B)$
which contradicts with \eqref{jol2.2} (since it implies that ${\rm
Re}(1 - \varphi_i(v_i)) > \alpha_0^2 /2$ for any $i\in I$).
\end{prf}

\medskip

We will show in the following that $A$ has property $(T)$
(respectively, strong property $(T)$) if and only if $U(A)$ has some
$(T)$-type property. Note that the argument of $(iii)\Rightarrow (i)$ in 
part (b) of the following result is adapted from that of Proposition 16 of
Chapter 1 of \cite{HV}.

\medskip

\begin{thm}
\label{U(A)-T} 
Let $\mathcal{SB}(U(A)) := \{(\mu, H) \in {\rm
Rep}(U(A)): (\mu,H)\leq (\pi_K,K)\ {\rm for\ some}\ K\in \sbm(A)\}$
(for ${\rm Rep}(U(A))$, we regard $U(A)$ as a discrete group).

\smnoind (a)  $A$ has property $(T)$ if and only if there exists
$(F,\epsilon) \in \F(U(A))\times \RP$ such that for any $K\in
\sbm(A)$ with $V_{\pi_K}(F,\epsilon)\neq \emptyset$, we have
$K^{U(A)}\neq \{0\}$.

\smnoind (b) The following statements are equivalent.

\begin{enumerate}
\item [\rm i.] $A$ has strong property $(T)$
\item [\rm ii.] There exists $(F,\epsilon) \in \F(U(A))\times
\RP$ such that for any $(\mu, H)\in \mathcal{SB}(U(A))$
and $\xi\in V_\mu(F,\epsilon)$, we have $P_H^{U(A)}(\xi) \neq 0$.
\item [\rm iii.] There exists $(F,\epsilon) \in \F(U(A))
\times \RP$ such that for any $(\mu, H)\in
\mathcal{SB}(U(A))$ with $V_\mu(F,\epsilon)\neq \emptyset$,
we have $H^{U(A)}\neq \{0\}$.
\end{enumerate}
\end{thm}
\begin{prf}
(a) This follows from Lemma \ref{F(U)}.

\smnoind (b) $(i) \Rightarrow (ii)$. 
Let $(F,\epsilon)\in
\F(U(A))\times \RP$ be a strong Kazhdan's pair for $(A, A, 1/2)$ (see Lemma \ref{F(U)}). 
Suppose that
$(\mu, H)\in \mathcal{SB}(U(A))$ and $K\in \sbm(A)$ such that $(\mu,
H) \leq (\pi_K, K)$. 
If $\xi \in V_\mu(F,\epsilon) \subseteq V_K(F,\epsilon)$, then 
$$\left\|\xi - P^A_K(\xi)\right\|\ <\ 1/2.$$ 
As $K = H\oplus H^\bot$ and both $H$ and $H^\bot$ are
invariant under $\pi_K$, we have $K^A = K^{U(A)} = H^{U(A)}\oplus
(H^\bot)^{U(A)}$. 
As $\xi\in H$, we know that $P^{U(A)}_{H^\bot}(\xi) =
0$ and $P^{U(A)}_H(\xi) = P^A_K(\xi) \neq 0$.

\smnoind $(ii) \Rightarrow (iii)$. This implication is clear.

\smnoind $(iii) \Rightarrow (i)$. 
Let $(F,\epsilon)$ be the pair
satisfying the hypothesis. 
For any $2\geq \alpha > 0$, we take
$\beta := \frac{\alpha\epsilon}{2}$. 
Suppose that $K\in \sbm(A)$, $H
= (K^A)^\bot = (K^{U(A)})^\bot$ and $\mu = \pi_K\!\mid\!_{H}$. 
Then
$(\mu, H) \in \mathcal{SB}(U(A))$. 
For any $\xi\in V_K(F, \beta)$,
we have $\xi = P^A_K(\xi) + \eta$ where $\eta\in H$. 
Assume that
$\epsilon \|\eta\| > \beta$ and put $\zeta := \eta/ \|\eta\|$. As
$$\|\mu(v)\zeta - \zeta\|\ =\ \|\pi_K(v)\xi - \xi\|/\|\eta\|\
<\ \beta/\|\eta\|\ <\ \epsilon \qquad (v\in F),$$
we have $\zeta\in V_\mu(F, \epsilon)$. Hence by the hypothesis,
$H$ contains a non-zero $\mu$-invariant vector which contradicts the
definition of $H$. Therefore, we must have $\|\eta\| \leq
\beta/\epsilon$ and hence $\left\|\xi - P^A_K(\xi)\right\| <
\alpha$.
\end{prf}

\medskip

\begin{rem}
(a) Let $\mathcal{SB}_0(U(A)) = \{(\mu, H) \in {\rm
Rep}(U(A)): (\mu,H)\leq (\pi_\FH,\FH)\}$. Then using the same
argument, one can show that a similar statement as Theorem
\ref{U(A)-T}(b) holds when $\mathcal{SB}(U(A))$ is replaced by
$\mathcal{SB}_0(U(A))$.

\smnoind
(b) Note that in the proof for $(i)\Rightarrow(ii)$ in Theorem \ref{U(A)-T}(b), we only need the existence of a strong Kazhdan's pair for $(A,A,1/2)$. 
\end{rem}

\smnoind Chi-Wai Leung, Department of Mathematics, The Chinese University of Hong Kong, Hong
Kong.

\smnoind \emph{Email address:} cwleung@math.cuhk.edu.hk

\medskip\noindent
Chi-Keung Ng, Chern Institute of Mathematics and LPMC, Nankai University, Tianjin 300071,
China.

\smnoind \emph{Email address:} ckng@nankai.edu.cn; ckngmath@hotmail.com

\bigskip


\begin{thebibliography}{99}
\bibitem{Ana}
C. Anantharaman-Delaroche, On Connes' property $(T)$ for von Neumann
algebras, Math. Japon. 32 (1987), 337-355.

\bibitem{Bek-unit}
M.B. Bekka, Kazhdan's property $(T)$ for the unitary group of a
separable Hilbert space, Geom. Funct. Anal. 13 (2003), 509-520.

\bibitem{Bek-T}
M.B. Bekka, Property $(T)$ for $C^*$-algebras, Bull. London Math. Soc.
38 (2006) 857-867.

\bibitem{BV} M.B. Bekka and A. Valette, Kazhdan's property $(T)$ and amenable representations, Math. Z. 212 (1993), 293-299.

\bibitem{Bro}
N.P. Brown, Kazhdan's property $(T)$ and $C^*$-algebras, Journal of
Functional Analysis 240 (2006) 290-296.

\bibitem{Con} 
A. Connes, Classification des facteurs, \emph{Operator algebras and applications, Part 2 (Kingston, Ont., 1980)}, Proc. Sympos. Pure Math. 38 (1982), 43-109. 

\bibitem{CJ}
A. Connes and V. Jones, Property $(T)$ for von Neumann algebras,
Bull. Lond. Math. Soc. 17 (1985), 57-62.

\bibitem{HV}
P. de la Harpe and A. Valette, La propri$\acute{\mbox{e}}$t$\acute{\mbox{e}}$ $(T)$ de
Kazhdan pour les groupes localement compacts, Ast\'{e}risque 175 (1989).

\bibitem{Jol-disc}
P. Jolissaint, Property $T$ for discrete groups in terms of their regular representation, Math. Ann. 297 (1993), 539-551.

\bibitem{Jol}
P. Jolissaint, On Property (T) for Pairs of Topological Groups,
Enseign. Math. (2) 51 (2005), 31-45.

\bibitem{Kaz} D. Kazhdan, Connection of the dual space of a group with the structure of its closed subgroups, Funct. Anal. Appl. 1 (1967), 63-65. 

\bibitem{Mar} G. Margulis, Finitely-additive invariant measures on Euclidean spaces, Erg. Th. and Dyn. Sys. 2 (1982), 383-396. 

\bibitem{Mur}
G. Murphy, \emph{$C\sp *$-algebras and operator theory}, Academic
Press, Boston MA, (1990).

\bibitem{Pat}
A.L.T. Paterson, Nuclear $C\sp *$-algebras have amenable unitary
groups, Proc. Amer. Math. Soc. 114 (1992), 719-721.

\bibitem{Popa}
S. Popa, On a class of type ${\rm II}\sb 1$ factors with Betti numbers invariants, Ann. of Math. 163 (2006), 809-899. 

\bibitem{V}
A. Valette, Old and new about Kazhdan's property $(T)$, Pitman research notes in Mathematics
series. Longman. vol.311 (1994), 271-333.
\end{thebibliography}
\end{document}